\documentclass[12pt,a4paper]{article}

\usepackage[left=2cm,right=1.5cm, top=2cm,bottom=3cm,bindingoffset=0cm]{geometry}
\usepackage[utf8]{inputenc}
\usepackage[T2A]{fontenc}
\usepackage[russian]{babel}
\usepackage{graphics}
\usepackage{amsfonts}
\usepackage{amsmath}
\usepackage{amssymb}
\usepackage{dsfont}
\usepackage{amsthm}
\usepackage{textcomp}
\usepackage{graphicx}
\usepackage{euscript}
\usepackage{mathrsfs}
\newtheorem{assert}{Assertion}
\newtheorem{lemma}{Lemma}
\newtheorem{consequence}{Corollary}
\newtheorem{defin}{Definition}
\newtheorem{theorem}{Theorem}
\author{Anastasia V. Parusnikova, Andrey V. Vasilyev}
\title{On Divergence of Puiseux Series Asymptotic Expansions of Solutions to the Third Painlev\'{e} Equation}%\thanks{This study was carried out within The National Research University Higher School of Economics Academic Fund Program in 2013-2014, research grant No. 12-01-0030.}}
\date{}

\begin{document}
\maketitle

\renewcommand{\abstractname}{Abstract} 
\begin{abstract}
    % In \cite{Ufa}, for formal asymptotic expansions of solutions of the third Painlev\'{e} equation into Puiseux series, it is shown that these series asymptotically approximate solutions of this equation in sectors with vertices at infinity. We give an upper estimate for the opening angle of the sector in which solutions can be approximated by  Puiseux series. Also, we give coefficients upper estimates, whose existence follows from the general theory of Gevrey formal series summability. 
    
    %В работе \cite{Ufa} для формальных асимптотических разложений решений третьего уравнения Пенлеве в ряды Пюизо показано, что эти ряды асимптотически приближают решения данного уравнения в секторах с вершиной в бесконечности. Для угла раствора сектора, в котором решения могут быть приближены рядами Пюизо дана оценка сверху, также приведены следующие из общей теории суммируемости по Жевре формальных рядов оценки сверху для коэффициентов.     
     
In this paper we present a family of values of the parameters of the third Painlev\'{e} equation such that Puiseux series formally satisfying this equation -- 
considered as series of $z^{2/3}$ -- are series of exact 
Gevrey order one. We prove the divergence of these series and provide analytic functions which are approximated by them %these power series
in sectors with the vertices at infinity. 

%Надо улучшать, а первый абзац думать, как перенести во вводный параграф


\textbf{MSC classes:} 	34m25, 34m55.

\end{abstract}

\section{Introduction}
~~~~Consider the third Painlev\'{e} equation
\begin{equation*}
\label{1}
w'' = \dfrac{(w')^2}{w} - \dfrac{w'}{z} + \dfrac{\alpha w^2 + \beta}{z} + \gamma w^3 + \dfrac{\delta}{w},
\end{equation*}
with $\gamma=0$, $\alpha, \beta, \delta \in \mathbb{C}$, $\alpha \delta \ne 0$. By making the following transform
\begin{equation*}
w(z) = z^{1 / 3}u(x), \: x = \frac{3}{2}z^{2 / 3},
\end{equation*}
we obtain the equation:
\begin{equation}
\label{P3_mod}
u''_{xx} = \dfrac{(u'_x)^2}{u} - \dfrac{u'_x}{x} +  \alpha u^2 + \dfrac{3}{2} \dfrac{\beta}{x} +  \dfrac{\delta}{u}.
\end{equation}

\smallskip
All formal power series satisfying equation (\ref{P3_mod}) near infinity have the %following
form
\begin{equation}
\label{series}
\sum_{n=0}^\infty a_n x^{-n},
\end{equation} where coefficients $a_n$ are determined by the formulae:
\begin{equation}
\label{coef}
a_0 = (-\delta / \alpha)^{1/3}; \quad a_1 = -\dfrac{\beta}{2\alpha a_0}; \quad a_2 = 0; \quad a_3 = -\dfrac{\beta}{6\alpha ^2 a_0^2} \biggl(\dfrac{\beta^2}{4 \delta} + 1\biggr);$$ $$a_n=\dfrac{1}{3\alpha a_0^2}\left(-\alpha\sum_{k=1}^{n-1}(a_k\sum_{j=0}^{n-k}a_ja_{n-k-j})-\alpha a_0\sum_{k=1}^{n-1}a_k a_{n-k}-\dfrac{3}{2}\beta a_{n-1}+\right.$$ $$ \left.+\sum_{k=1}^{n-2}(2k^2+k(2-n))a_k a_{n-2-k}\right), \ n\geqslant 4.
\end{equation} 

\begin{defin}
The series $\sum\limits_{n=0}^\infty a_n z^{-n}$ is called a series of Gevrey order $k$~$\cite{Balser}$ if there exist constants $k, M, C>0$ such that
\begin{equation}
\label{estimate}
|a_n|\leqslant C (n!)^{1/k}M^n~\forall n \in \mathbb{N}.
\end{equation}
\end{defin}
 
\begin{defin}
The series $\sum\limits_{n=0}^\infty a_n z^{-n}$ is called a series of  exact Gevrey  order $k$~$\cite{Ramis}$ if it is of Gevrey order $k$ and there exists no $k'>k$ such that it is of Gevrey order $k'$.
\end{defin} 
 
In \cite{Ufa} it is shown that series~(\ref{series}) is a series of Gevrey order one and as is proved in the book~\cite{Gromak}, the series~(\ref{series}) is a rational function %-- that is a series $\sum\limits_{n=0}^\infty a_n \left(\frac{3}{2}\right)^nz^{(2n+1)/3}$ with the coefficients $a_n$~(\ref{coef}) is an algebraic function -- 
iff $\beta\neq 0, \delta = -\beta^2/(4k)^2, k \in \mathbb{Z}\setminus\{0\}$ or $\beta =0$: for such values of parameters the estimates (\ref{estimate}) are not precise and can be improved.
The aim of the present paper is to show that there also exist the values of the parameters of the third Painlev\'{e} equation for which the series (\ref{series}) is of exact Gevrey order one, %, we show that there also exist the values of the parameters of the equation $P_3$ for that %the lower estimate of the coefficients $a_n$ is of the similar form, meaning 
%there exist constants $M_1, C_1, \ell>0$ such that $|a_n|\geq C_1(n-\ell)!M_1^n~\forall n \in \mathbb{N}$ 
hence the series (\ref{series}) diverges.  

In the second section we find such values of parameters and prove the following
\begin{theorem}
Series {\upshape (\ref{series})} with coefficients~{\upshape (\ref{coef})} with parameters $\alpha, \beta \in \mathbb{C}$, $\alpha\beta\neq 0$, $\delta =-\beta^2/2$
is of exact Gevrey order one.
\end{theorem}
In the third section we construct analytic functions being approximated by series (\ref{series}) for parameters of the third Painlev\'{e} equation $\alpha=-1/32, \ \beta=-1/4, \ \delta =-1/32$.

\section{Proof of Theorem 1}
~~~~Here we prove the divergence of power series~(\ref{series}) with coefficients~(\ref{coef}) considering the fixed values of the parameters of equation~(\ref{P3_mod}) $\alpha=-1/32, \ \beta=-1/4, \ \delta =-1/32$. Assume that the branch of the cube root is fixed so that $a_0 = -1$ with the parameters given. Note that under the above conditions all the coefficients $a_n \in \mathbb{R}$.

%Let us
Write out the first coefficients of the series needed for the further calculations:
$$a_0=-1, a_1=4, a_2=0, a_3=\dfrac{64}{3}, a_4=\dfrac{256}{3},a_5=2048.$$
With $n \geqslant 6$ the coefficients have the form
\begin{equation}
\notag
\label{coef_param}
a_n=-\dfrac{1}{3}\sum_{k=1}^{n-1}(a_k\sum_{j=0}^{n-k}a_ja_{n-k-j})+\dfrac{1}{3} \sum_{k=1}^{n-1}a_k a_{n-k}-4 a_{n-1}-\dfrac{32}{3}\sum_{k=1}^{n-2}(2k^2+k(2-n))a_k a_{n-2-k};
\end{equation}
we rearrange them assuming that $n \geqslant 9$ and considering the values of the first coefficients calculated earlier. We obtain
\begin{equation}
\label{coef_simplify}
a_n=4 a_{n-1}-16 a_{n-2}-\dfrac{1}{3}\sum_{k=3}^{n-6}\left(a_k\sum_{j=3}^{n-k-3}a_j a_{n-k-j}\right)+\sum_{k=3}^{n-3}a_k a_{n-k}-4\sum_{k=3}^{n-4}a_k a_{n-1-k}-$$ $$-\dfrac{32}{3}\sum_{k=3}^{n-5}(2k^2+k(2-n))a_k a_{n-2-k}-\dfrac{128}{3}(n-4)^2a_{n-3}+\dfrac{32}{3}(n-2)^2a_{n-2}.
\end{equation}

\begin{assert}
For the coefficients $a_n$  $(\ref{coef_simplify})$ of series $(\ref{series})$, the following holds:
\begin{equation}
\begin{aligned}
&(A)~~ a_n \leqslant n! \left(\dfrac{32}{3}\right)^{\frac{n}{2}-1},~n\geq 4;\\
&(B)~~a_n \geqslant \dfrac{32}{3}(n-3)(n-2)a_{n-2},~n\geq 5;\\
&(C)~~a_n \geqslant \sum_{k=3}^{n-2}a_k a_{n-k+1},~n\geq 5;\\
&(D)~~a_n >0,~n\geq 3.
\end{aligned}
\end{equation}
\end{assert}

\textbf{Proof of Assertion 1} (induction proof).
%Induction base: 
Inequality (A) holds with $n=4, \dots, 8$, inequalities (B) and (C) hold with $n=5, \dots, 8$, inequality (D) holds with $n=3, \dots, 8$.

%Induction step: 
Let $n \geqslant 9$ and
\begin{equation}
\begin{aligned}
\label{step}
(A)~~&a_m \leqslant m! \left(\dfrac{32}{3}\right)^{\frac{m}{2}-1},~m=4, \ldots, n-1;\\
(B)~~&a_m \geqslant \dfrac{32}{3}(m-3)(m-2)a_{m-2},~m=5, \ldots, n-1;\\
(C)~~&a_m \geqslant \sum_{k=3}^{m-2}a_k a_{m-k+1},~m=5, \ldots, n-1;\\
(D)~~&a_m >0,~m=3, \ldots, n-1.
\end{aligned}
\end{equation}
We prove that
$a_n \leqslant n! \left(\dfrac{32}{3}\right)^{\frac{n}{2}-1};~~a_n \geqslant \dfrac{32}{3}(n-3)(n-2)a_{n-2};~~a_n \geqslant \sum\limits_{k=3}^{n-2}a_k a_{n-k+1};~~a_n >0.$

\begin{lemma}
\begin{equation}
\label{lemma1}
\sum_{k=3}^{n-5}(2k^2+k(2-n))a_k a_{n-2-k}\geqslant 0, \: n \geqslant 20.
\end{equation}
\end{lemma}
Proof of Lemma 1.
 Split the monomials on the left-hand side (LHS) of inequality (\ref{lemma1}) into pairs and sum up the monomial corresponding to index $k$ and the monomial corresponding to index $n-2-k$, where $k<(n-2)/2$. If $n$ is odd, then the LHS sum of inequality (\ref{lemma1}) splits into the above pairs; if $n$ is even, then, besides the above pairs, there remains the expression $a^2_{(n-2)/2}(2((n-2)/2)^2-(n-2)/2\cdot(2-n))=0$; that is, in order to prove inequality (\ref{lemma1}), it is sufficient to check non-negativity of the coefficient of $a_k a_{n-2-k}$, since, by induction hypothesis and base, $a_k>0, a_{n-2-k}>0$. This coefficient is equal to
$$2 k^2+k(2-n)+2(n-2-k)^2+(n-2-k)(2-n)=(2k-(n-2))^2\geqslant 0,$$
which completes the proof of Lemma 1. \qed

\begin{lemma}
\begin{equation}
\notag
\label{lemma2}
	k!(n-k)!\leqslant 4!(n-4)! \mbox{ with } k =4, \dots, n-4.
\end{equation}
\end{lemma}
Proof of Lemma 2.
$C_n^k\geq C_n^4, k =4, \dots, n-4$ that is $$\dfrac{n!}{k!(n-k)!}\geqslant \dfrac{n!}{4!(n-4)!}, k =4, \dots, n-4;$$
therefore $\dfrac{1}{k!(n-k)!}\geqslant \dfrac{1}{4!(n-4)!}, \: k =4, \dots, n-4$, whence the statement of Lemma 2 follows. \qed

\textbf{Proof of inequality (A) of Assertion 1}. Replacing terms on the right-hand side (RHS) of equality (\ref{coef_simplify}) with equal or greater expressions, and using statements of Lemmas 1 and 2 and induction hypothesis (A), we obtain 
$$a_n\leqslant4a_{n-1}+\sum_{k=3}^{n-3}a_k a_{n-k}+\dfrac{32}{3}(n-2)^2a_{n-2}=4a_{n-1}+\dfrac{128}{3}a_{n-3}+\sum_{k=4}^{n-4}a_k a_{n-k}+\dfrac{32}{3}(n-2)^2a_{n-2}\leqslant$$
$$\!\leqslant\!4{\left(\dfrac{32}{3}\right)\!}^{\frac{n-3}{2}}\!(n-1)!+\dfrac{128}{3}{\left(\dfrac{32}{3}\right)\!}^{\frac{n-5}{2}}\!(n-3)!+{\left(\dfrac{32}{3}\right)\!}^{\frac{n-4}{2}}\!\sum_{k=4}^{n-4}k!(n-k)!+\dfrac{32}{3}\left(\dfrac{32}{3}\right)^{\frac{n-4}{2}}\!(n-2)^2(n-2)!\leqslant$$
$$\leqslant\left(\dfrac{32}{3}\right)^{\frac{n-4}{2}}\left(4\sqrt{\dfrac{32}{3}}(n-1)!+4\sqrt{\dfrac{32}{3}}(n-3)!+4!(n-4)!(n-7)+\dfrac{32}{3}(n-2)^2(n-2)!\right)=\tilde{A}.$$

Thus, to prove inequality (A) it is sufficient to check that
\begin{equation}
\label{A+}
\tilde{A}\leqslant n!\left(\dfrac{32}{3}\right)^{\frac{n}{2}-1}.
\end{equation}

Dividing the LHS and RHS of inequality (\ref{A+}) by the positive $8(n-4)!\left(\dfrac{32}{3}\right)^{\frac{n}{2}-2}$, and rearranging the terms, we see that to complete the proof of inequality (A) of Assertion 1 it remains to show that
$$3(n-7)\leqslant \dfrac{4}{3}(n-3)\left(n(n-1)(n-2)-\sqrt{\dfrac{3}{2}}(n-2)(n-1)-\sqrt{\dfrac{3}{2}}-(n-2)^3\right)=$$ $$=\dfrac{4}{3}(n-3)\left(\left(3-\sqrt{\dfrac{3}{2}}\right)n^2-\left(10-3\sqrt{\dfrac{3}{2}}\right)n+8-3\sqrt{\dfrac{3}{2}}\right).$$ 

The last inequality holds, since for the function $$f(x)=\left(3-\sqrt{\dfrac{3}{2}}\right)x^3-\left(19-6\sqrt{\dfrac{3}{2}}\right)x^2+\left(\dfrac{143}{4}-12\sqrt{\dfrac{3}{2}}\right)x+9\sqrt{\dfrac{3}{2}}-\dfrac{33}{4},$$ it is true that $f(x)>0$ with $x>0$, therefore, $f(n)>0$ with $n\geq 9$, which finishes the proof of inequality~(A).

\textbf{Proof of inequality (B) of Assertion 1}. First we prove the following 
%several auxiliary 
inequalities:

~~1) $a_{n-1}\geqslant \dfrac{32}{3}(n-4)(n-3)a_{n-3} \geqslant \dfrac{32}{3}(n-4)^2a_{n-3}$ with $n \geqslant 4$ by induction hypothesis (B) with $m=n-1$  and by inequality (D) with $m=n-3$;
$$\mbox{2)}~\sum\limits_{k=3}^{n-4}a_ka_{n-k-1}\leqslant  a_{n-2} \mbox{~by induction hypothesis (C) with~} m=n-2;\hspace{2.65cm}$$
$$\mbox{3)}~\sum_{k=3}^{n-6}\left(a_k\sum_{j=3}^{n-k-3}a_j a_{n-k-j}\right)\stackrel{\mbox{\footnotesize with }m=5, \dots, n-4 \mbox{\footnotesize ~in inequality (C) }}{\leqslant}\hspace{4.85cm}$$
$$
\leqslant \sum_{k=3}^{n-6}a_k a_{n-1-k}\stackrel{\mbox{\footnotesize by inequality (D) }}{\leqslant} \sum_{k=3}^{n-4}a_k a_{n-1-k} \stackrel{\mbox{\footnotesize with }m= n-2 \mbox{\footnotesize ~in inequality (C) }}{\leqslant}  a_{n-2};
$$
\begin{equation}
\label{ner_4}
\mbox{4)}~\sum_{k=3}^{n-3}a_k a_{n-k}\geqslant\dfrac{32}{3}\sum_{k=3}^{n-5}(2k^2+k(2-n))a_k a_{n-2-k}.\hspace{5.7cm}
\end{equation}

To prove (\ref{ner_4}), we decrease its LHS by subtracting from it two positive (by induction hypothesis (D) (\ref{step})) monomials with $k=n-4$, $k=n-3$; then we use inequality (B) (\ref{step}) for every $a_{n-k}$ in its LHS by taking $m=5, \dots, n-3$, that is
$$\sum_{k=3}^{n-3}a_k a_{n-k}\geqslant \sum_{k=3}^{n-5}a_k a_{n-k}\geqslant \dfrac{32}{3}\sum_{k=3}^{n-5}(n-k-3)(n-k-2)a_k a_{n-2-k}.$$

Thus, to prove (\ref{ner_4}) it is sufficient to check that
$$\dfrac{32}{3}\sum_{k=3}^{n-5}(n-k-3)(n-k-2)a_k a_{n-2-k}\geqslant \dfrac{32}{3}\sum_{k=3}^{n-5}(2k^2+k(2-n))a_k a_{n-2-k}, \mbox{ or}$$
\begin{equation}
\label{coefB}
 \sum_{k=3}^{n-5}\left((n-k-3)(n-k-2)-2k^2-2k+kn\right)a_k a_{n-2-k}\geqslant 0.
\end{equation} 

Split the monomials in the LHS of (\ref{coefB}) into pairs and sum up the monomial corresponding to index $k$ and the monomial corresponding to index $n-2-k$, where $k<(n-2)/2$. Suppose $n$  odd, then the LHS sum of inequality (\ref{coefB}) splits into the above pairs. Suppose $n$  even, then, besides the above pairs, there remains the expression $\dfrac{a^2_{(n-2)/2}}{4}\left(n-2\right)(n-4)\geqslant 0$ with $n\geqslant 4$, that is in order to prove inequality~(\ref{ner_4}), it is sufficient to check non-negativity of the coefficient of $a_k a_{n-2-k}$, $3\leq k < (n-2)/2$ \: since, by induction hypothesis, $a_k>0, a_{n-2-k}>0$. This coefficient is equal to
$$(n-k-3)(n-k-2)-2k^2-2k+kn+k(k-1)-2(n-2-k)^2-2(n-2-k)+(n-2-k)n=$$ $$=2(k-1)(n-k-3)+n-4 \geqslant 0$$ 
for the indices $k$ considered in (\ref{coefB}).

Due to inequalities 1) -- 4) and formula (\ref{coef_simplify}), the coefficient
\begin{equation}
\label{coef_B}
\begin{aligned}
a_n&\geqslant\dfrac{128}{3}(n-4)^2a_{n-3}-16 a_{n-2}-\dfrac{1}{3}a_{n-2}-4a_{n-2}-\dfrac{128}{3}(n-4)^2a_{n-3}+\dfrac{32}{3}(n-2)^2a_{n-2}=\\
&=\left(\dfrac{32}{3}(n-2)^2-\dfrac{61}{3}\right)a_{n-2}\geqslant\dfrac{32}{3}(n-3)(n-2) a_{n-2} \mbox{ with } n \geqslant 4,
\end{aligned}
\end{equation}
which proves inequality (B) of Assertion 1.

\begin{lemma}
\begin{equation}
\label{lowerbound}
a_n \geqslant 4\left(\dfrac{32}{3}\right)^{\frac{n}{2}-1}(n-2)! \mbox{\; with \:} n\geqslant 5.
\end{equation}
\end{lemma}
Proof of Lemma 3. Assume $n$ is even. Then we use $\left(\frac{n}{2} - 2\right)$ times with $m = 6,8, \dots, n$ the already proven inequality (B) of Assertion 1:
\begin{gather*}
a_n \geqslant \left(\dfrac{32}{3}\right)^{\frac{n}{2} - 2}(n-2)(n-3)\dots (n - (n-3)) a_4
= 4\left(\dfrac{32}{3}\right)^{\frac{n}{2}-1}(n-2)! \mbox{\; with \:} n\geqslant 6.
\end{gather*}
If $n$ is odd, we use $\frac{n-3}{2}$ times inequality (B) of Assertion~1  with $m = 5, 7, \dots, n$:
\begin{equation*}
a_{n}\geqslant{\left(\dfrac{32}{3}\right)\!}^{\frac{n-3}{2}}(n-2)(n-3)\dots(n-(n-2))a_3
 \!=2\left(\dfrac{32}{3}\right)^{\frac{n-1}{2}}(n-2)! \mbox{\; with \:} n\geqslant 5.
\end{equation*}

From $2\biggl(\dfrac{32}{3}\biggr)^{\frac{n-1}{2}}(n-2)! > 4\biggl(\dfrac{32}{3}\biggr)^{\frac{n}{2}-1}(n-2)! \mbox{\; with \:} n\geqslant 5$, we obtain that inequality~(\ref{lowerbound}) holds with all $n \geqslant 5$. Lemma 3 is proved. \qed

%\textbf{To prove inequality (C) of Assertion 1}, let us use inequality (\ref{coef_B}). We have\\ $a_n\geq \left(\dfrac{32}{3}(n-2)^2-\dfrac{61}{3}\right)a_{n-2}$, that is why it is sufficient to check that
\textbf{To prove inequality (C) of Assertion 1}, we use inequality (\ref{coef_B}). As we have\\ $a_n\geqslant \left(\dfrac{32}{3}(n-2)^2-\dfrac{61}{3}\right)a_{n-2}$,  it is sufficient to check that
$$\left(\dfrac{32}{3}(n-2)^2-\dfrac{61}{3}\right)a_{n-2}\geqslant \sum_{k=3}^{n-2}a_ka_{n-k+1}=\sum_{k=4}^{n-3}a_ka_{n-k+1}+2 a_3a_{n-2}, \mbox{ that is}$$
\begin{equation}
\label{coef_C}
\dfrac{32}{3}\left(n^2-4n-\dfrac{61}{32}\right)a_{n-2}\geqslant\sum_{k=4}^{n-3}a_ka_{n-k+1}.
\end{equation}

For upper estimate %bound?
of the RHS of inequality (\ref{coef_C}) we use induction hypothesis (A) for each $a_\ell$, and also use an analogue of Lemma 2: 
\begin{equation}
\label{coef_C_1}
\sum_{k=4}^{n-3}a_k a_{n-k+1}\leqslant \sum_{k=4}^{n-3}\left(\dfrac{32}{3}\right)^{\frac{n-3}{2}}k!(n-k+1)!=\left(\dfrac{32}{3}\right)^{\frac{n-3}{2}}\left(\sum_{k=5}^{n-4}k!(n-k+1)!+2 \cdot 4!(n-3)!\right)\leqslant$$ 
$$\leqslant\left(\dfrac{32}{3}\right)^{\frac{n-3}{2}}\biggl(5!(n-4)!(n-8)+2 \cdot 4!(n-3)!\biggr).
\end{equation}

For lower estimate %bound?
of the LHS of  inequality (\ref{coef_C}) we use Lemma 3 and obtain 
\begin{gather}
\label{coef_C_2}
\dfrac{32}{3}\left(n^2\!-\!4n\!-\!\dfrac{61}{32}\right)a_{n-2}\geqslant 4{\left(\dfrac{32}{3}\right)\!}^{\frac{n}{2}-1}\left(n^2\!-\!4n\!-\!\dfrac{61}{32}\right)(n-4)! \mbox{\; with \:} n\geqslant 7.
\end{gather}

Due to above inequalities (\ref{coef_C_1}), (\ref{coef_C_2}) to check inequality (\ref{coef_C}) it is sufficient to prove that the following inequality holds with $n \geqslant 9$:
\begin{equation}
\label{coef_C_next}
4{\left(\dfrac{32}{3}\right)\!}^{\frac{n}{2}-1}\left(n^2\!-\!4n\!-\!\dfrac{61}{32}\right)(n-4)!\geqslant \left(\dfrac{32}{3}\right)^{\frac{n-3}{2}}\biggl(5!(n-4)!(n-8)+2 \cdot 4!(n-3)!\biggr).
\end{equation}

Divide both LHS and RHS of inequality (\ref{coef_C_next}) by $(n-4)!\left(\dfrac{32}{3}\right)^{\frac{n-3}{2}}$, carry all the terms to the LHS; we prove non-negativity of the expression
\begin{gather*}
4 \sqrt{\dfrac{32}{3}}\left(n^2\!-\!4n\!-\!\dfrac{61}{32}\right) - 120(n-8) - 48(n-3) =\\ 
= 16 \sqrt{\dfrac{2}{3}}n^2 - \biggl(168 + 64\sqrt{\dfrac{2}{3}}\biggr)n + 1104 -\dfrac{61}{2}\sqrt{\dfrac{2}{3}} >0 \mbox{ with } n\geqslant 9.
\end{gather*}

	\textbf{Inequality (D) of Assertion 1} follows from inequality (C) and induction hypothesis (D), since 
$$ a_n\geqslant \dfrac{32}{3}(n-3)(n-2)a_{n-2}>0 \mbox{\; with \:} n\geqslant 5.$$
This completes the proof of Assertion 1. \qed

\begin{consequence}
For the coefficients $a_n$ {\upshape (\ref{coef})} of series {\upshape (\ref{series})} with the fixed values of the parameters $\alpha=-1/32, \ \beta=-1/4, \ \delta =-1/32$ the following holds:
\begin{equation}
\label{estimates}
4\left(\dfrac{32}{3}\right)^{\frac{n}{2}-1}(n-2)! \ \leqslant a_n \ \leqslant  \left(\dfrac{32}{3}\right)^{\frac{n}{2}-1}n! \mbox{\; with \:} n\geqslant 5.
\end{equation}
\end{consequence}

%%\smallskip
%%\includegraphics{graph.jpg}

%%{\itshape In order to understand how precise estimates~{\upshape  (\ref{estimates})} are, we consider the functions \\$\left(32/3\right)^{(x/2) - 1}\Gamma(x+1)$ and $4\left(32/3\right)^{(x/2) - 1}\Gamma(x-1)$, whose values coincide with estimates~{\upshape  (\ref{estimates})} with $x \in \bbN, \ x\geq5$, and plot the logarithms of these functions together with the logarithms of $a_n$ values.}

\begin{lemma}
If another branch of the cube root is fixed in the calculation of coefficients {\upshape (\ref{coef})} of series {\upshape (\ref{series})}, that is $\tilde a_0 = \omega  a_0$, where $\omega^3 = 1$, then
\begin{equation}
\label{roots}
\notag
\tilde a_n = \omega^{n+1}a_n \mbox{\; with \:} n\geq 0.
\end{equation}
\end{lemma}

The proof is by induction on $n$. \qed

\smallskip
Lemma 4 together with estimates~(\ref{estimates}) immediately imply 
%This lemma shows that all the three formal solutions diverge.

%Доказательство. Докажем по индукции по $n$. База индукции: при $n = 0, 1, 2$ и 3 утверждение~(\ref{roots}) выполнено. 

%Шаг индукции: пусть $n\geq 4$ и $\tilde a_m = \omega^{m+1}a_m$ при $m = 0, \dots, m-1$. Тогда 
%\begin{gather*}
%\tilde a_n=\dfrac{1}{3\alpha \tilde a_0^2}\left(-\alpha\sum_{k=1}^{n-1}(\tilde a_k\sum_{j=0}^{n-k}\tilde a_j \tilde a_{n-k-j})-\alpha \tilde a_0\sum_{k=1}^{n-1}\tilde a_k \tilde a_{n-k}-\dfrac{3}{2}\beta \tilde a_{n-1}+\right.\\
%\left.+\sum_{k=1}^{n-2}(2k^2+k(2-n))\tilde a_k \tilde a_{n-2-k}\right) = \dfrac{\omega}{3\alpha a_0^2}\left(-
%\alpha \omega^{n+3} \sum_{k=1}^{n-1}(a_k\sum_{j=0}^{n-k}a_j a_{n-k-j})\right.-\\
%-\left.\alpha \omega^{n+3} a_0\sum_{k=1}^{n-1}a_k  a_{n-k}-\dfrac{3}{2} \omega^{n} \beta \tilde a_{n-1} + \omega^{n} \sum_{k=1}^{n-2}(2k^2+k(2-n))a_k a_{n-2-k}\right) = \omega^{n+1}a_n.
%\end{gather*}

\begin{consequence}
Series {\upshape (\ref{series})} with coefficients~{\upshape (\ref{coef})} with the fixed values of the parameters\\ $\alpha=-1/32, \ \beta=-1/4, \ \delta =-1/32$
 diverges for any branch of the cube root.
\end{consequence}

\begin{consequence}
Series {\upshape (\ref{series})} with coefficients~{\upshape (\ref{coef})} with parameters $\alpha, \beta \in \mathbb{C}$, $\alpha\beta\neq 0$, $\delta =-\beta^2/2$
diverges.
\end{consequence}

The proof can be easily obtained from the following assertion~\cite{Gromak}:

Let $w = \varphi(z)$ be a solution of the third Painlev\'{e} equation with given values $\alpha, \beta, \gamma, \delta$ of the parameters, then the function $\sigma_1 \varphi(\sigma_2 z)$
where $\sigma_1$, $\sigma_2\in \mathbb{C}, \: \sigma_1\sigma_2\neq 0$ is a solution of the third Painlev\'{e} equation with the parameters  $\alpha \sigma_1^{-1} \sigma_2,$ $\beta \sigma_1 \sigma_2,$ $\gamma\sigma_1^{-2} \sigma_2^{2},$ $\delta\sigma_1^{2} \sigma_2^{2}.$ \qed

%Пусть $w = \varphi(z)$ -- решение $P_3$ при заданных значениях $\alpha, \beta, \gamma, \delta$ параметров уравнения, тогда 
%\begin{equation}
%\label{mapGromak}
%\sigma_1 \varphi(\sigma_2 z),
%\end{equation}
%где $\sigma_1$, $\sigma_2\in \mathbb{C}, \sigma_1\sigma_2\neq 0$ является решением $P_3$ при параметрах уравнения  $\alpha \sigma_1^{-1} \sigma_2,$ $\beta \sigma_1 \sigma_2,$ $\gamma\sigma_1^{-2} \sigma_2^{2},$ $\delta\sigma_1^{2} \sigma_2^{2}.$

%Ряд {\upshape (\ref{series})} с коэффициентами {\upshape (\ref{coef})} расходится одновременно с решением уравнения $P_3$, которое получается из этого ряда при обратной замене переменных. Таким образом, ряд~{\upshape (\ref{series})} с коэффициентами~{\upshape (\ref{coef})} будет расходиться при любых параметрах $\alpha, \beta$ и $\delta$, которые можно перевести в \ $(-1/32), (-1/4)$ и $(-1/32)$ соответственно при помощи преобразования~(\ref{mapGromak}).

%Возьмем произвольные $\alpha, \beta \in \mathbb{C}$, $\alpha\beta\neq 0$ и подберем такие $\sigma_1$ и $\sigma_2$, что $\sigma_1^{-1} \sigma_2 = -1/(32\alpha)$, $\sigma_1 \sigma_2 = -1/(4\beta)$. Тогда при $\delta =-\beta^2/2$ преобразование~(\ref{mapGromak}) переводит параметры $\alpha, \beta$ и $\delta$ в \ $(-1/32), (-1/4)$ и $(-1/32)$ соответственно, что завершает доказательство.

\section{On Borel and Laplace transforms}
~~~~As is proved in \cite{Ufa} solutions to the third Painlev\'{e} equation considered as functions of $z^{2/3}$ are asymptotically approximated of Gevrey order one by the series (\ref{series}) in the sectors with the vertices at infinity with opening not larger then $\pi$. As we see from Corollary 3 series (\ref{series}) with parameters of the equation $\alpha, \beta \in \mathbb{C}$, $\alpha\beta\neq 0$, $\delta =-\beta^2/2$ diverges and does not present an analytic solution to the third Painlev\'{e} equation. We construct an analytic functions being approximated of Gevrey order one by series (\ref{series}) in the same sectors as mentioned above obtaining the Borel sum of the series (\ref{series}) and then applying formal Laplace transform to it. %Here we briefly recall this construction. 

Firstly, consider series 
\begin{equation}
\label{series_b_n}
\hat{g}(x)=\sum\limits_{n=0}^\infty b_n x^{-n},
\end{equation}
with $b_n=\dfrac{a_n}{n!}$, this series is a formal Borel transform of series (\ref{series}). Series (\ref{series_b_n}) converges for $|x|>R$. We calculate $R$ using the estimates from Corollary 1 and applying Cauchy-Hadamard theorem:
$$\frac{1}{R}=\overline{\lim_{n \rightarrow \infty}}\sqrt[n]{|b_n|}=\sqrt{\dfrac{32}{3}}.$$

Then we apply Borel-Ritt-Gevrey theorem \cite{Balser}: given an arbitrary sector $S$ of opening at most $\pi$ with $d \in \mathbb{R}$ being the bisecting direction of $S$ and $\rho \in \mathbb{R}, |\rho|>R$ the formal finite Laplace transform to series (\ref{series_b_n}) 
$$g\left(\frac{1}{x}\right)=\dfrac{1}{x}\int\limits_{0}^{ e^{id}/\rho}\left(\sum_{n=0}^\infty b_n u^n e^{-u/x}\right)du,$$
where integrating is along $\textrm{Arg}\,u = d$, gives us an analytic function $g(x)$ asymptotically approximated by series~(\ref{series}) in the sector $S\cap \{x:|x|>R\}$. %Хочется понять, нужно ли здесь пересечение или асимптотика во всём секторе, как и говорится в теореме из книги Бальзера.

Hence, we conclude that the difference between $g(x)$ and a solution to the third Painlev\'{e} equation in the given sector $S$ is a series of Gevrey order one. %(not necessarily of exact order one).

%Our next purpose is to investigate the analytic functions obtained and to study the nonlinear Stokes phenomena concerned with solutions 

\renewcommand{\refname}{References}

\vskip 1 cm
\textbf{Affiliations.}\\

Anastasia V. Parusnikova\\
National Research University Higher School of Economics,\\
%Moscow Institute of Electronics and Mathematics,\\
34 Tallinskaya str., 123458, Moscow, Russia\\
e-mails: parus-a@mail.ru, aparusnikova@hse.ru

\vskip 1 cm
Andrey V. Vasilyev\\
National Research University Higher School of Economics,\\
6 Usacheva str., 119048, Moscow, Russia\\
e-mail: vasiljev.andr@gmail.com
\end{document}